\documentclass[a4paper,11pt,twoside,reqno]{article}
\RequirePackage[latin1]{inputenc}
\RequirePackage[T1]{fontenc}
\RequirePackage[english]{babel}
\RequirePackage{amsmath}
\RequirePackage{amssymb}
\numberwithin{equation}{section}
\newcommand{\N}{\mathbf{N}}
\renewcommand{\phi}{\varphi}

\newenvironment{proof}[1][Proof :]{\textbf{#1} }{\hfill $\blacksquare$}

\def\hs{\hspace}
\def\O{\Omega}

\def\eps{\varepsilon}
\def\a{\alpha}
\def\tl{\widetilde}
\def\disp{\displaystyle}

\def\E{{\bf E}}
\def\Prb{{\bf P}}
\def\qs{\forall\;}

\def\td{\longrightarrow}
\def\cal{\mathcal}
\newtheorem{theo}{Theorem}[section]
\newtheorem{lem1}{Lemma}[section]

\title{\bf Besov  regularity  of  the  uniform  empirical  process}
\author{Gane Samb Lo \and Ahmadou Bamba SOW}
\date{\today}
\begin{document}
\maketitle
\noindent{\footnotesize LERSTAD, UFR S.A.T, Universit\'e Gaston
Berger, BP 234, Saint-Louis, SENEGAL. \\ Phone: (221) 33 961 23 40.
Email : {\tt ganesamblo@yahoo.com,  ahbsow@gmail.com}}
%LERSTAD, UFR S.A.T, Universit\'e Gaston Berger, B.P. 234
%Saint-Louis, S\' en\'egal. \\{\sl email: [gslo,
%abnabsow]@yahoo.fr}\vspace{1cm}
%
\begin{abstract}
The paths of Brownian motion have been widely studied in the recent
years relatively in Besov spaces $B_{p, \infty}^\a$. The results are
the same as to the Brownian bridge. In fact these regularities
properties are established in some
sequence spaces $S_{p, \infty}^\a$ using an isomorphisim between them and $B_{p, \infty}^\a$.\\
In this note, we are concerned with the regularity of the paths of
the continuous version of the uniform empirical process in the space
$S_{p, \infty}^\a$ and in one of his separable sub space $S_{p,
\infty}^{\a, 0}$ for  a  suitable  choice   of  $\a$   and $p$.
\end{abstract}
AMS Classification subject: 1991, 60G17, 60G15, 60F15. 
\section{Introduction}
Let   $U_1, U_2, ..., U_n, \dots  $ be  a  sequence  of   i.i.d
$\mathcal{U}(0, 1)$ random variables. For a fix integer $n \ge 1$
  we consider  the empirical  distribution   function $\tl F_n$ of  the  sample  $U_1, U_2, \dots, U_n$  defined
  by
$$
\qs  0  \le  s  \le 1, \quad   \tl F_n(s) =  \frac{1}{n} \,
\sum_{i=1}^n 1_{(-\infty, \, s ]} (U_i)
$$
and     for   $j \ge 0, \; k = 1, ..., 2^j$ the triangular sequence
\begin{equation}
\tl\a_{jk}^n = 2^{j/2} \bigg[2 \,\tl\a_n (\frac{k - 1/2}{2^j})  -
\tl\a_n (\frac{k - 1}{2^j})  - \tl\a_n (\frac{k}{2^j} )\bigg]
\label{alpha-n}
\end{equation}
where   $\tl \a_n $ is  the associated empirical process defined
   by  $\tl\a_n(s) = \sqrt n (\tl
F_n (s) - s ), \quad 0\le s\le 1$. \\
Our   motivation    in  the study  of this sequence is
  given  by  previous   works    on   the   regularity   of  the  paths   of  the  Brownian  motion
 in  Besov  spaces   $B_{p,
\infty}^\a $ given    by
$$
B_{p, \infty}^\a =  \{f  \in  L^p ([0, 1]) :   \sup_t \frac{w_p(f,
t)}{t^\a}  <  \infty\}
$$
where   for  any  $1 \le  p < \infty$,
$$ w_p(f, t) =  \sup_{|h|\le t} \bigg(\int_{I_h} |f(x  - h)  -  f(x)|^p\, d x \bigg)^{1/p}
;  \quad  I_h  =  \{ x  \in [0, 1],  \;   x  -  h  \in  [0, 1]\}.
$$
The  space  $B_{p, \infty}^\a$ endowed   with  the  following norm
$$
||f ||_{p, \infty}^\a  =  \sup\{ |f_0| ,  |f_{1}|,  \;  \sup_j
2^{-j(1/2 - \a - 1/p)} \bigg( \sum_{k=1}^{2^j}
|f_{jk}|^p\bigg)^{1/p} \}
$$
where
$$f_0  =  f(0),  \;  f_1  =  f(1)  - f(0), \; \quad  f_{jk} = 2^{j/2}
\bigg[2 f (\frac{k - 1/2}{2^j})  - f (\frac{k - 1}{2^j}) - f
(\frac{k}{2^j} )\bigg] $$
%and   the space
is    a    Banach  space.

It  is  well  known    thanks   to Ciesielski et al \cite{CRK}, that
there   exists   an  isomorphism    between   such     spaces  and
the  Banach  spaces   of  sequences  $(S_{p,
\infty}^\a,||\cdot||_{p, \infty}^\a)$   defined by
%The study of the regularity of stochastic process in
% such  spaces  is  improve    by  Ciesielski's et al \cite{CRK} isomorphism betwennL'\'etude des processus stochastiques relativement \`a ces espaces
%est rendue accesssible gr\^ace \`a l'isomorphisme construite par
%Ciesielski et al \cite{CRK}   entre ces espaces et l'espace de
%suites $S_{p, \infty}^\a$  d\'efini par
%
$$
\{\mu = (\mu_{jk},  \;  j \ge 0, \;  k = 1,  ..., 2^j) /  \; \,
||\mu ||_{p, \infty}^\a \,  <  \,  \infty\}
%\sup_j 2^{-j(1/2 - \a - 1/p)} \bigg(\sum_{k=1}^{2^j} |\mu_{jk}|^p\bigg)^{1/p}  <  \infty\}
$$
%muni de la norme
%
where  $$ ||\mu ||_{p, \infty}^\a  =  \sup\{ |\mu_0| ,  |\mu_{1}|,
\;  \sup_j 2^{-j(1/2 - \a - 1/p)} \bigg( \sum_{k=1}^{2^j}
|\mu_{jk}|^p\bigg)^{1/p} \}.
$$
Their subsets   $B_{p, \infty}^{\a, 0}$ (respectively  $S_{p,
\infty}^{\a, 0}$) of functions   $  f  \in B_{p, \infty}^\a$ (resp
of sequences $(\mu_{jk}) \in S_{p, \infty}^\a $) such  that $w_p(f,
t) = \circ(t^\a)$  as   $t \td 0 $ (resp $ 2^{-j(1/2 - \a - 1/p)}
\bigg( \sum_{k=1}^{2^j} |\mu_{jk}|^p\bigg)^{1/p} \td  0 $  as $j \td
\infty$)
are   separable   Banach spaces. \\
Thanks    to  this  isomorphism, Roynette \cite{ROY}  proved  that
for  $p \ge 2$ and  $\a < 1/2$, the   Brownian   path   $(W_t, 0\le
t\le 1)$ belongs    almost  surely  in   $B_{p, \infty}^\a$ but not
   in  $ B_{p, \infty}^{\a, 0}$ by  establishing
$$
 \sup_j \bigg( 2^{-j}  \sum_{k=1}^{2^j} |g_{jk}|^p \bigg)^{1/p}\; <  \infty \; \; \mbox{and}\;\;
 \liminf_{j\td +\infty} \bigg(2^{-j} \sum_{k=1}^{2^j} |g_{jk}|^p \bigg)^{1/p}  \; >  0
$$
where   for  all    $ j \ge 0$   and   $k = 1,  ..., 2^j, \quad
g_{jk} =
2^{j/2} \bigg[2 W (\frac{k - 1/2}{2^j})  - W (\frac{k - 1}{2^j})  - W (\frac{k}{2^j} )\bigg] $. \\
This  result  can  be  extended     to  the   Brownian bridge $b_t =
W_t - t W_1, \; t \in [0, 1]$ which   is   closely  related   to the
  uniform empirical  process. Moreover   Koml\'os  et  al \cite{KMT} show   that    on    a   suitable
probability space $(\O, {\cal A}, \Prb)$,  there   exists a
 sequence   of  i.i.d. $\mathcal{U}(0, 1)$ random variables $ U_1, U_2, ..., $
 and   a   sequence  of   Brownian bridges
$\{b_n(t),   0\le  t \le 1\}$ such  that almost  surely
$$
\disp \limsup_{n \td +\infty} \; \; \frac{\sqrt n}{\log n}\;\;  \sup_{0\le t \le 1} |\a_n(t)  -  b_n(t)|  <  \infty.
$$
Our   aim  is  to  investigate   Roynette's result   for the
Brownian
 bridge   to  the  continuous  version  of   the  empirical process. We   successfully   get  a result   for $1\le p \le 2 $ and   $\a =  1/2$.
%
%
%\section{Results}
\section{Empirical process  and Besov Spaces}
In   order   to   face  to  the   lack  of  smoothness   of  the
classical  empirical   process,      we   first  recall the
following result :
\begin{lem1}
For  every  $n\ge 1$, the   empirical  distibution    process $\tl
\a_n$ admits a continuous version $\a_n$.
\end{lem1}
\begin{proof}
It is  well    known    one  can  express   the    distribution
empirical    function  $\tl  F_n(\cdot)$  in  terms  of  the order
statistics  $U_1^{(n)}   \le  U_2^{(n)}  \le   \dots  \le U_n^{(n)}$
of  the  sample  $U_1, U_2, \dots, U_n$   as   follows
\begin{eqnarray*}
\tl  F_n (s)  =
\begin{cases}
0,  \quad \; U_1^{(n)}  >  s,  \\
\disp  \frac{k}{n},   \quad   U_k^{(n)}   \le   s   <   U_{k
+1}^{(n)}, \quad k
= 1, 2,  \dots, n   - 1,  \\
1,  \quad   U_n^{(n)}   \le    s.
\end{cases}
\end{eqnarray*}
Let   us   consider the   function   $F_n(\cdot) $   defined  for
every $0 \le  s  \le  1$ by
$$  F_n (s) = \tl F_n(U_k^{(n)}) +   2   \Bigg( s      -
\frac{U_k^{(n)}   + U_{k+1}^{(n)}}{2} \Bigg)\Bigg(\frac{\tl
F_n(U_{k+1}^{(n)})    - \tl F_n(U_k^{(n)})}{U_{k+2}^{(n)}  -
U_{k}^{(n)}}\Bigg),  $$ if
 $$ \frac{U_k^{(n)}   +  U_{k+1}^{(n)}}{2}   \le  s   \le
\frac{U_{k+1}^{(n)}   +  U_{k+2}^{(n)}}{2} .$$ It  is   easy   to
see that for every  $n  \ge 1$, \begin{equation}
\label{estimate1}\sup_{0\le s\le 1} |F_n (s) - \tl F_n(s) | \; \le
\frac{1}{n}.\end{equation} As a consequence of   \eqref{estimate1},
we deduce that  $F_n$   is  a   continuous    version  of   $\tl
F_n$   and   the  process
  $\a_n(s) = \sqrt n (F_n (s) - s ), \;  0\le s\le 1$   is    a
continuous version of the associated empirical process $\tl\a_n(s) =
\sqrt n (\tl F_n (s) - s ), \; 0\le s\le 1$.
% where   $F_n^*$ stands for  the  continuous  version  of  $F_n$.
\end{proof}\\ \\
We  are   now  in  position   to   formulate   our  main  results
\begin{theo}
For  every   $n  \ge 1$, the    process   $\a_n$   satisfy   almost
surely
$$
{\{\a_{jk}^n\}}_{j,k} \in S_{2, \infty}^{1/2} \quad \mbox{and }\quad
  {\{\a_{jk}^n\}}_{j,k} \notin S_{2, \infty}^{1/2, 0}.
$$
\end{theo}
\noindent\begin{proof} Let  us    consider   the   triangular
sequence  given   by   \eqref{alpha-n} (replacing  $\tl \a_n$ by
 $\a_n$), we   deduce  thanks   to the distribution empirical process
 that
$$
\qs j \ge 0,  \;  \qs k = 1, ..., 2^j,  \;\;  \; \a_{jk}^n =  \frac{2^{j/2}}{\sqrt n} \sum_{i=1}^n  Z_{jk}(i)
$$
where
$$
\qs  i  = 1, ..., n, \;  \;\;\;  Z_{jk}(i) = Z_{jk}(U_i) =  1_{[\frac{k - 1}{2^j}, \frac{k - 1/2}{2^j} [} (U_i)  -  1_{[\frac{k - 1/2}{2^j}, \frac{k}{2^j}[} (U_i)
.$$
Notice  that     for  any     $i = 1, ..., n,  \;\; Z_{jk}(i) \in \{
1, 0 , - 1 \} $   respectively    with   probability  $\disp
\frac{k}{2^j}, 0, \frac{k}{2^j}$. We  deduce   that    for   any $i
= 1, ..., n, \;\; Z_{jk}(i)$  is   centered   random    variable
 with   variance    $2^{-j} $. \\
Let  us   define   $\disp G_{jk} =  |\a_{jk}^n|^2  =  \frac{2^j}{n} H_{jk} $
where    $\disp H_{jk}  =  \bigg(\sum_{i=1}^n  Z_{jk}(i)\bigg)^2$. \\
Using the   fact   that   for  any  fix   $k $  and   $ i \neq h$  the  random  variables  $ Z_{jk}(i) $  and    $ Z_{jk}(h)  $
are  independent we   deduce  that   %
$$
\E (H_{jk}) =  \E  \bigg(\sum_{i=1}^n  Z^2_{jk}(i)   +  \,
\sum_{i\neq h}^n  Z_{jk}(i)  Z_{jk}(h) \bigg )  =  \frac{n}{2^j}.
 $$
which  imlplies  in  particular   $\E (G_{jk})  = 1 $. \\
Futhermore  for  any  $j \ge 0$   and   $k = 1, ..., 2^j $, we have
\begin{eqnarray*}
H_{jk}^2  &=&  \bigg(\sum_{i=1}^n  Z_{jk}^2(i) \bigg)^2  + 2
\sum_{l=1}^n  Z_{jk}^2(l)\sum_{i\neq h}^n  Z_{jk}(i)  Z_{jk}(h) +
\bigg(\sum_{i\neq h}^n  Z_{jk}(i)  Z_{jk}(h) \bigg)^2 \\
&=& \sum_{i=1}^n  Z_{jk}^4(i)  + 2 \sum_{i < h}^n  Z_{jk}^2(i)
Z_{jk}^2(h)   +  2 \sum_{i\neq h}^n  Z_{jk}^2(i)  Z_{jk}^2(h)
%2 \sum_{l=1}^n \sum_{i\neq h}^n Z_{jk}^2(l)  Z_{jk}(i)  Z_{jk}(h) \\
+ A_{jk}^{(1)}   + A_{jk}^{(2)}
%2 \sum_{l\neq m}^n \sum_{i\neq h}^n   Z_{jk}(i)  Z_{jk}(h)  Z_{jk}(l)  Z_{jk}(m)
 \end{eqnarray*}
where for   every $j \ge 0$   and   $k = 1, ..., 2^j $,
$$
\disp A_{jk}^{(1)} = 2 \sum_{l=1}^n \sum_{i\neq h}^n Z_{jk}^2(l)
Z_{jk}(i) Z_{jk}(h) \; \; \mbox{and}\;\;    A_{jk}^{(2)} =   2
\sum_{l\neq m}^n \sum_{i\neq h}^n   Z_{jk}(i)  Z_{jk}(h)  Z_{jk}(l)
Z_{jk}(m)
$$
It is  easy  to see   that for   every $j \ge 0$   and   $k = 1,
..., 2^j, \, A_{jk}^{(1)}  $  is   a   centered   random  variable
and   $A_{jk}^{(2)} $   satisfies
\begin{eqnarray*}
A_{jk}^{(2)}&=& 4   \sum_{i < h}^n \bigg[\sum_{l < m}^n   Z_{jk}(i)  Z_{jk}(h)  Z_{jk}(l)  Z_{jk}(m) \bigg] \\
&=&   4   \sum_{i < h}^n \bigg[\sum_{(l,m)=(i, h) }   Z_{jk}(i)  Z_{jk}(h)  Z_{jk}(l)  Z_{jk}(m) \bigg] + 4 \sum_{i < h}^n \bigg[\sum_{(l,m)\neq(i, h) }   Z_{jk}(i)  Z_{jk}(h)  Z_{jk}(l)  Z_{jk}(m) \bigg]
\end{eqnarray*}
The   expectation   of  the   last  sum   vanish thanks   to  the
independence   of   the   random   variables. Hence  there   exists
a   constant  $c > 0$
 which   may  change   from  line   to  line such  that   $\E (A_{jk}^{(2)})   =\disp   c \; \E \sum_{i < h}^n Z^2_{jk}(i)\,
Z^2_{jk}(h)$. We    deduce   that for   every $j \ge 0$ and $k = 1,
..., 2^j$,
%\begin{align*}
$$\E (H_{jk}^2)   =    \sum_{i =1}^n   \E (Z_{jk}^4(i))  + c \; \E
\sum_{i < h}^n   Z_{jk}^2(i) \,  Z_{jk}^2(h)  =  \frac{n}{2^j} (1 +
c  \frac{n - 1}{2^j}) $$ which  implies  in particular for every $j
\ge 0$   and   $k = 1, ..., 2^j$,
$$ Var (G_{jk}) =  \frac{2^j}{n}[1  + \frac{n(c - 1)   -  c}{2^j}].
 $$
Elsewhere  we  have
\begin{eqnarray*}
Var (\sum_{k=1}^{2^j} G_{jk})  &=& \sum_{k=1}^{2^j} \, Var (G_{jk})  + 2 \sum_{1=k < k' \le 2^j}\frac{2^{2j}}{n^2} \, cov (H_{jk}\,H_{jk'} ) \\
H_{jk}\,H_{jk'}  & =&   \bigg( \sum_{i =1}^n  Z_{jk}(i) Z_{jk'}(i)    + \sum_{i \neq h }^n  Z_{jk}(i) Z_{jk'}(h)  \bigg)^2
\end{eqnarray*}
Notice  that   for  $k  \neq k'$, the  product   $Z_{jk}(i)
Z_{jk'}(i)$ is null. This   implies for every $j \ge 0$   and   $k
\ne k' \in  \{1, ..., 2^j\}$,
\begin{eqnarray}
 H_{jk}\,H_{jk'}  & =&   \bigg( \sum_{i \neq  h}^n  Z_{jk}(i) Z_{jk'}(h) \bigg)^2
= \sum_{i \neq  h}^n\sum_{l \neq  m}^n   Z_{jk}(i) Z_{jk'}(h) Z_{jk}(l) Z_{jk'}(m) \label{card}\\
&=& \sum_{card \{i, h\}\cap \{l, m\} = 2} \; Z_{jk}(i) Z_{jk'}(h) Z_{jk}(l) Z_{jk'}(m) \nonumber\\
&+& \sum_{card \{i, h\}\cap \{l, m\} <  2} \; Z_{jk}(i) Z_{jk'}(h) Z_{jk}(l) Z_{jk'}(m)\nonumber
\end{eqnarray}
So   two   cases  can  be   investigated : \\
if  $card \{i, h\}\cap \{l, m\} <  2$, extracting  one  random  variable   $Z_{jk}(i)$, the  expectation   of  the last term in \eqref{card} is  null. \\
if   $card \{i, h\}\cap \{l, m\} =  2$, we   have   either  ($l = i$
and  $m = h$) or   ($l = h$ and   $m = i$).   In  this  last  case
the   product   is   equal  to  zero. It  remains
\begin{eqnarray*}
\E(H_{jk}\,H_{jk'})  &=&   \sum_{i \neq  h}^n \E[Z_{jk}^2(i) Z_{jk'}^2(h)]  = \frac{n(n - 1)}{2^{2j}} \\
Var (\sum_{k=1}^{2^j} G_{jk})  &=& \sum_{k=1}^{2^j} \,\frac{2^j}{n}(1  + \frac{3n  - 4}{2^j}) \;  +\;   2 \sum_{1=k < k' \le 2^j}\frac{2^{2j}}{n^2} \,\bigg( \frac{n(n - 1)}{2^{2j}}    -  (\frac{n}{2^j})^2\bigg) \\
&=&  2^{2j} \eps_{nj},   \;\;\;\;\mbox{where}\;\;\;\;  \eps_{nj} =
\frac{1}{2^j} (3  - \frac{3}{n})
\end{eqnarray*}
Exploiting  Bienaymé-Tch\'ebychev  inequality,     we  obtain the
  following   estimate for  every  $n  \in \N$ and
$$
\qs  j \ge 0,  \;\;\; \Prb( |2^{-j}  \, \sum_{k=1}^{2^j} |\a_{jk}^n|^2  \;  -  1  |  \; \ge \frac{1}{2}) \;\le  4 \; \eps_{nj}
$$
Therefore   thanks   to  Borel-Cantelli   lemma, we deduce that for
any $n  \in \N$
$$
 \frac{1}{2}    \le \;  2^{-j}  \, \sum_{k=1}^{2^j} |\a_{jk}^n|^2
\;\, \le  \frac{3}{2}\;\;\;  p.s,\;\;  j \; \mbox{large enough.}
$$
Hence
$$
\sup_j \;  \bigg( 2^{-j}  \, \sum_{k=1}^{2^j} |\a_{jk}^n|^2
\bigg)^{1/2} \; <  \infty \hs{1cm}\mbox{and}\hs{1cm} \liminf_{j \td
+ \infty} \bigg( 2^{-j}  \, \sum_{k=1}^{2^j} |\a_{jk}^n|^2
\bigg)^{1/2} \; >  0 \; \; a.s.
$$
\end{proof} \\
{\bf Remark:}
\begin{enumerate} \item Theses   results
show that for any $n$, the sequence  $(\a_{jk}^n), \; j\ge 0, \; \;
k = 1, ..., 2^j $ belongs
    in  the   space  $S_{2, \infty}^{1/2}$  and  not   in   $S_{2, \infty}^{1/2, 0}$ a.s.
\item The   first   result  of   our   theorem  can  be  extended  to  $1\le p \le 2$  since   the   $L^p$  norm is   incerasing  in  $p$.
\end{enumerate}


\begin{thebibliography}{5}
\bibitem{Boufoussi} B. Boufousssi (1994). {\em Espaces de Besov: Caract\'erisations et applications},
Th\`ese de l'Universit\'e Henri Poincar\'e, Nancy I, France.
\bibitem{CRK} Z. Ciesielski, B. Roynette, G. Kerkyacharian (1993).
 {\em Quelques espaces fonctionnels associ\'es \`a des processus gaussiens}, Studia Mathematica, 107, p. 171-204.
\bibitem{KMT} J. Koml\'os, M. Major, G. Tusnàdy (1975). {\em Weak convergence and
embedding}. In colloquia Math.Soc. Janos.Boylai.Limit Theorems of
Probability Theory, 149-165. Amsterdam, North-Holland.
%\bibitem{reference} trouver refer sur  net.  {\em Strong approximations  of empirical   processes  by  Gaussian Processes}.

\bibitem{ROY} B. Roynette (1993). {\em Mouvement brownien et espaces de Besov}, Stochastics and Stochastics Reports, 43, p. 221-260.
\end{thebibliography}
 \end{document}